\documentstyle[12pt]{article}
\def\mineappendix{
        \setcounter{section}{1}
        \setcounter{subsection}{0}
        \def\thesection{\Alph{section}}
        \def\sectionap{\@startsection  {section}{1}{\z@}
                        {-3.5ex plus-1ex minus-.2ex} {0ex plus.2ex}
                        {\reset@font\Large\bf  Appendix:  \, }
                        }
        }
\makeatother
\def\Proclaim #1. #2\par{\bigbreak\noindent{\sc#1.\enspace}{\it#2}\par}

\font\Bbbfont=msbm10
\newfam\msbfam
\textfont\msbfam=\Bbbfont
\def\Bbb#1{{\fam\msbfam\relax#1}}

\newcommand{\eqref}[1]{equation~(\ref{#1})}

\newcommand{\gwii}[1]{\left< \hspace{-2pt} \left< \, #1 \,
        \right>  \hspace{-2pt} \right>_{0}}

\newcommand{\gwiione}[1]{\left< \hspace{-2pt} \left< \, #1 \,
        \right> \hspace{-2pt} \right>_{1}}

\newcommand{\gwiitwo}[1]{\left< \hspace{-2pt} \left< \, #1 \,
        \right> \hspace{-2pt} \right>_{2}}

\newcommand{\gwig}[1]{\left< \, #1 \, \right>_{g}}
\newcommand{\gwiig}[1]{\left< \hspace{-2pt} \left< \, #1 \,
    \right> \hspace{-2pt} \right>_{g}}

\newcommand{\grav}[2]{\tau_{#1}(\gamma_{#2})}
\newcommand{\grava}[1]{\tau_{#1}(\gamma_{\alpha})}

\newcommand{\gravb}[1]{\tau_{#1}(\gamma_{\beta})}

\newcommand{\ga}{\gamma_{\alpha}}
\newcommand{\gua}{\gamma^{\alpha}}
\newcommand{\gb}{\gamma_{\beta}}
\newcommand{\gub}{\gamma^{\beta}}
\newcommand{\gm}{\gamma_{\mu}}
\newcommand{\gum}{\gamma^{\mu}}

\newcommand{\vs}{{\cal S}}
\newcommand{\vx}{{\cal X}}
\newcommand{\vd}{{\cal D}}
\newcommand{\ve}{{\cal E}}
\newcommand{\vf}{{\cal F}}
\newcommand{\vw}{{\cal W}}
\newcommand{\vv}{{\cal V}}
\newcommand{\vu}{{\cal U}}
\newcommand{\vl}{{\cal L}}
\newcommand{\vg}{{\cal G}}

\newcommand{\bvs}{{\, \overline{\cal S}\,}}
\newcommand{\bvx}{{\, \overline{\cal X}\,}}

\newcommand{\bvw}{{\, \overline{\cal W} \,}}

\newcommand{\bvl}{{\, \overline{\cal L} \,}}

\newcommand{\qp} {\circ}

\newtheorem{lem}{Lemma}[section]
\newtheorem{cor}[lem]{Corollary}

\newtheorem{pro}[lem]{Proposition}

\title{Idempotents on the big phase space}
\author{Xiaobo Liu \thanks{Research partially supported by
            Alfred P. Sloan Research Fellowship and National Science
        Foundation research grant}}
\date{}

\begin{document}
\maketitle

Let $M$ be a compact symplectic manifold.
In Gromov-Witten theory, the space $H^{*}(M; {\Bbb C})$ is called the
{\it small phase space}. The so called
large quantum cohomology provides a ring structure
on each tangent space of the small phase space. Together with the intersection
pairing, this defines a Frobenius manifold structure on the small phase space.
As a vector space, the small phase space has a natural flat coordinate system.
However, this coordinate system is hard to work with in practice.
When the quantum cohomology is generically semisimple, there is another
coordinate system on the small phase space, called the
canonical coordinate system,
whose coordinate vector fields are idempotents of the quantum product.
This coordinate system played a vital role in solving the genus-1
Gromov-Witten potential in \cite{DZ1}.  It seems hard to use the
canonical coordinate system to study Gromov-Witten invariants
of genus bigger than 1. The main reason for this is that starting
from genus 2, all known universal equations for Gromov-Witten invariants
involve descendant classes, and therefore are not
equations on the small phase space (cf. \cite{G2} and \cite{BP}).
Instead, they are partial differential equations on the
{\it big phase space}, which can be think of as a product
of infinitely many copies of the small phase space.
The main purpose of this paper is to set up a frame work on the big phase space
which play the role of the canonical coordinate system on the small phase
space.

Quantum product of vector fields on the big phase space can be
naturally defined using the genus-0 Gromov-Witten potential (cf.
\cite{L2}). This product coincides with the large quantum product
when restricted to the small phase space. But on the big phase
space, the quantum product does not have identity. This implies
that we can not extend the canonical coordinate system to the big
phase space in a natural way. Instead, we can still consider
vector fields on the big phase space which are idempotents of the
quantum product. When the quantum product is generically
semisimple, idempotents span the space of primary vector fields
on the big phase space. We can apply an operator $T$ defined in
\cite{L2} repeatedly to idempotents and obtain a nice frame for
the tangent bundle of the big phase space. This frame is not
commutative with respect to the Lie bracket, and therefore can
not be coordinate vector fields for any coordinate system. But
they also have many good properties which greatly simplify known
universal equations. In this paper, we will study basic properties
of idempotents on the big phase space and indicate how to apply
them to study Gromov-Witten invariants of genus 1 and 2.
They can also be applied to study higher genus
Gromov-Witten invariants once the corresponding universal
equations are found. Since universal equations also appear in the
Gromov-Witten theory of spin curves, the frame work developed
here can also be applied to that setting.

Part of the work in this paper was done when the author visited
IPAM at Los Angeles, MSRI at Berkeley, and IHES in France. The
author would like to thank these institutes for hospitality.

\section{Idempotents}
\label{sec:idem}

Let's first recall the definition
of the quantum product on the big phase space as given in \cite{L2}.
Let $M$ be a compact symplectic manifold. For simplicity, we assume
$H^{\rm odd}(M; {\Bbb C}) = 0$.
The {\it big phase space} is by definition the product of infinite copies of
$H^{*}(M; {\Bbb C})$, i.e.
\[  P := \prod_{n=0}^{\infty} H^{*}(M; {\Bbb C}). \]
Fix a basis $\{ \gamma_{1}, \ldots, \gamma_{N} \}$ of
$H^{*}(M; {\Bbb C})$ with $\gamma_{1} = 1$ being the identity of the ordinary
cohomology ring of $M$. Then we denote the corresponding basis for
the $n$-th copy of $H^{*}(M; {\Bbb C})$ in $P$ by
$\{\tau_{n}(\gamma_{1}), \ldots, \tau_{n}(\gamma_{N}) \}$.
We call $\grava{n}$ a {\it descendant} of $\gamma_{\alpha}$ with descendant
level $n$.
We can think of $P$ as an infinite dimensional vector space with basis
$\{ \grava{n} \mid 1 \leq \alpha \leq N, \, \, \, n \in {\Bbb Z}_{\geq 0} \}$
where ${\Bbb Z}_{\geq 0} = \{ n \in {\Bbb Z} \mid n \geq 0\}$.
Let
$(t_{n}^{\alpha} \mid 1 \leq \alpha \leq N, \, \, \, n \in {\Bbb Z}_{\geq 0})$
be the corresponding coordinate system on $P$.
For convenience, we identify $\grava{n}$ with the coordinate vector field
$\frac{\partial}{\partial t_{n}^{\alpha}}$ on $P$ for $n \geq 0$.
If $n<0$, $\grava{n}$ is understood as the $0$ vector field.
We also abbreviate $\grava{0}$ as $\gamma_{\alpha}$.
Any vector field of the form $\sum_{\alpha} f_{\alpha} \ga$, where $f_{\alpha}$
are functions on the big phase space, is called a {\it primary vector field}.
We use $\tau_{+}$ and $\tau_{-}$ to denote the operator which shift the level
of descendants, i.e.
\[ \tau_{\pm} \left(\sum_{n, \alpha} f_{n, \alpha} \grava{n}\right)
    = \sum_{n, \alpha} f_{n, \alpha} \grava{n \pm 1} \]
where $f_{n, \alpha}$ are functions on the big phase space.

We will use the following {\it conventions for notations}:  All
summations are over the entire meaningful ranges of the indices
unless otherwise indicated. Let
\[ \eta_{\alpha \beta} = \int_{M} \gamma_{\alpha} \cup
    \gamma_{\beta}
\]
 be the intersection form on $H^{*}(M, {\Bbb C})$.
We will use $\eta = (\eta_{\alpha \beta})$ and $\eta^{-1} =
(\eta^{\alpha \beta})$ to lower and raise indices.
For example,
\[ \gua :=  \eta^{\alpha \beta} \gb.\]
Here we are using the summation convention that repeated
indices (in this formula, $\beta$) should be summed
over their entire ranges.

Let
\[ \gwig{\grav{n_{1}}{\alpha_{1}} \, \grav{n_{2}}{\alpha_{2}} \,
    \ldots \, \grav{n_{k}}{\alpha_{k}}} \]
be the genus-$g$ descendant Gromov-Witten invariant associated
to $\gamma_{\alpha_{1}}, \ldots, \gamma_{\alpha_{k}}$ and nonnegative
integers $n_{1}, \ldots, n_{k}$
(cf. \cite{W}, \cite{RT}, \cite{LiT}).
The genus-$g$
generating function is defined to be
\[ F_{g} =  \sum_{k \geq 0} \frac{1}{k!}
         \sum_{ \begin{array}{c}
        {\scriptstyle \alpha_{1}, \ldots, \alpha_{k}} \\
                {\scriptstyle  n_{1}, \ldots, n_{k}}
                \end{array}}
                t^{\alpha_{1}}_{n_{1}} \cdots t^{\alpha_{k}}_{n_{k}}
    \gwig{\grav{n_{1}}{\alpha_{1}} \, \grav{n_{2}}{\alpha_{2}} \,
        \ldots \, \grav{n_{k}}{\alpha_{k}}}. \]
This function is understood as a formal power series of
$t_{n}^{\alpha}$.

Introduce
a $k$-tensor
 $\left< \left< \right. \right. \underbrace{\cdot \cdots \cdot}_{k}
        \left. \left. \right> \right> $
defined by
\[ \gwiig{{\cal W}_{1} {\cal W}_{2} \cdots {\cal W}_{k}} \, \,
         := \sum_{m_{1}, \alpha_{1}, \ldots, m_{k}, \alpha_{k}}
                f^{1}_{m_{1}, \alpha_{1}} \cdots f^{k}_{m_{k}, \alpha_{k}}
        \, \, \, \frac{\partial^{k}}{\partial t^{\alpha_{1}}_{m_{1}}
            \partial t^{\alpha_{2}}_{m_{k}} \cdots
            \partial t^{\alpha_{k}}_{m_{k}}} F_{g},
 \]
for vector fields ${\cal W}_{i} = \sum_{m, \alpha}
f^{i}_{m, \alpha} \, \frac{\partial}{\partial t_{m}^{\alpha}}$ where
$f^{i}_{m, \alpha}$ are  functions on the big phase space.
We can also view this tensor as the $k$-th covariant derivative
of $F_{g}$ with respect to the trivial connection on $P$. This tensor is called the {\it $k$-point
(correlation) function}.
For any vector fields $\vw_{1}$ and $\vw_{2}$ on the big phase space,
the {\it quantum product}
 of $\vw_{1}$ and $\vw_{2}$ is defined by
\[ \vw_{1} \qp \vw_{2} := \gwii{\vw_{1} \, \vw_{2} \, \gua} \ga. \]
This is a commutative and associative product. But it does not have
an identity. For any vector field $\vw$ and integer $k \geq 1$,
$\vw^{k}$ is understood as the $k$-th power of $\vw$ with respect to this
product.

Let
\[{\vx} := - \sum_{m, \alpha} \left(m + b_{\alpha} - b_{1} -1
                        \right)\tilde{t}^{\alpha}_{m} \, \grava{m}
        - \sum_{m, \alpha, \beta}
        {\cal C}_{\alpha}^{\beta}\tilde{t}^{\alpha}_{m} \,
        \gravb{m-1}
\]
be the {\it Euler vector field} on the big phase space $P$, where
$\tilde{t}^{\alpha}_{m} = t^{\alpha}_{m} - \delta_{m,1} \delta_{\alpha, 1}$,
\[ b_{\alpha} = \frac{1}{2}({\rm dimension \,\,\, of \,\,\,} \ga)
   - \frac{1}{4} ({\rm real \,\,\, dimension \,\,\, of \,\,\,} M)
   + \frac{1}{2}\]
 and the matrix ${\cal C} = ( {\cal C}_{\alpha}^{\beta})$
is defined by $ c_{1}(V) \cup \ga = {\cal C}_{\alpha}^{\beta} \,
\gb$. For smooth projective varieties, the dimension
of $\ga$ can be replaced by twice of the holomorphic dimension of $\ga$
in the definition of $b_{\alpha}$. Moreover, the basis $\{
\gamma_{1}, \ldots, \gamma_{N}\}$ of $H^{*}(V, {\Bbb C})$ can be
chosen in a way such that the following holds: If $\eta^{\alpha
\beta} \neq 0$ or $\eta_{\alpha \beta} \neq 0$, then $b_{\alpha} =
1-b_{\beta}$.

The quantum multiplication by $\vx$ is an endomorphism on the space
of primary vector fields on $P$.
If this endomorphism has distinct eigenvalues at generic points,
we call $P$ {\it semisimple}. In this case, let $\ve_{1}, \ldots, \ve_{N}$
be the eigenvectors with corresponding eigenvalues $u_{1}, \ldots, u_{N}$,
i.e.
\[ \vx \qp \ve_{i} = u_{i} \ve_{i} \]
for each $i = 1, \cdots, N$. If we consider $\ve_{i}$ as a vector field on
$P$, then $u_{i}$ can be considered as a function on $P$.
By the associativity of the quantum product
\[ \vx \qp (\ve_{i} \qp \ve_{j}) = (\vx \qp \ve_{i}) \qp \ve_{j}
    = u_{i} \, \ve_{i} \qp \ve_{j}. \]
By the commutativity of quantum product, we also have
\[ \vx \qp (\ve_{i} \qp \ve_{j}) = u_{j} \, \ve_{i} \qp \ve_{j}. \]
Therefore at generic points where $u_{i} \neq u_{j}$ for $i \neq j$,
\[ \ve_{i} \qp \ve_{j} = 0 \]
and $\ve_{i} \qp \ve_{i}$ must be proportional to $\ve_{i}$.

Let
\[ {\vs} := - \sum_{m, \alpha} \tilde{t}^{\alpha}_{m}
        \grav{m-1}{\alpha} \]
be the {\it string vector field} on $P$. We define
\[ \bvw = \vw \qp \vs \]
for any vector field $\vw$ on $P$. The string equation implies that
the vector field $\bvs$ is the identity for the quantum product restricted
to the space of primary vector fields.
Since $\{ \ve_{i} \mid i=1, \ldots, N\}$ form a basis of the space of primary
vector fields, $\bvs$ is a linear combination of these vector fields.
Therefore the fact $\bvs \qp \ve_{i} = \ve_{i}$ implies that
$\ve_{i} \qp \ve_{i} \neq 0$ because otherwise $\vs \qp \ve_{i} = 0$.
Therefore we can normalize $\ve_{i}$ so that it is an idempotent of the
quantum product. Hence we can assume that
\[ \ve_{i} \qp \ve_{j} = \delta_{ij} \ve_{i} \]
for all $i$ and $j$.
We have
\begin{equation} \label{eqn:bvs}
     \bvs = \sum_{i=1}^{N} \ve_{i}.
\end{equation}
and
\begin{equation} \label{eqn:bvx}
 \bvx^{k} = \sum_{i=1}^{N} u_{i}^{k} \, \ve_{i}.
\end{equation}
for $k \geq 1$. One can verify these equations by multiplying both
sides of the equations by $\ve_{j}$ for arbitrary $j$.

Let $\nabla$ be the covariant derivative on $P$ of the trivial
flat connection with respect to the standard coordinates $\{
t_{n}^{\alpha} \}$. The compatibility of the quantum product and
this flat connection is given by the formula
\begin{equation} \label{eqn:DerProd}
 \nabla_{\vv} (\vw \qp \vu)
= (\nabla_{\vv} \vw) \qp \vu
    + \vw \qp (\nabla_{\vv} \vu) + \gwii{\vv \, \vw \, \vu \, \gua} \ga
\end{equation}
for any vector fields $\cal U$, $\cal V$ and $\cal W$ (cf.
\cite[Equation (8)]{L2}).

\begin{lem} \label{lem:deridem}
For any vector field $\vw$,
\[ \nabla_{_{\vw}} \ve_{i}
= -2 \gwii{\vw \, \ve_{i} \, \ve_{i} \, \gua} \ga \qp \ve_{i}
        + \gwii{\vw \, \ve_{i} \, \ve_{i} \, \gua} \ga \]
for each $i$.
\end{lem}
{\bf Proof}: Since $\ve_{i} = \ve_{i} \qp \ve_{i}$, by
\eqref{eqn:DerProd},
\[ \nabla_{_{\vw}} \ve_{i} = \nabla_{_{\vw}} (\ve_{i} \qp \ve_{i})
    = 2 (\nabla_{\vw} \ve_{i}) \qp \ve_{i} +
        \gwii{\vw \, \ve_{i} \, \ve_{i} \gua} \ga. \]
Multiplying both sides by $\ve_{i}$, we have
\[ (\nabla_{_{\vw}} \ve_{i} ) \qp \ve_{i}
    = - \gwii{\vw \, \ve_{i} \, \ve_{i} \gua} \ga \qp \ve_{i}. \]
Plugging this formula into the above equation, we obtain the desired formula.
$\Box$

For any vector field $\vw$, define
\[ T(\vw) := \tau_{+}(\vw) - \vs \qp \tau_{+}(\vw). \]
The operator $T$ was introduced in \cite{L2} to simplify
topological recursion relations for Gromov-Witten invariants. It
corresponds to the $\psi$ classes in the relations in the
tautological ring of moduli space of stable curves. In some
sense, repeatedly applying $T$ to a vector field will trivialize
its action on genus-$g$ generating functions. Here are some basic
properties of $T$ (cf. \cite[Section 1]{L2}): For any vector fields $\vw_{i}$,
\begin{eqnarray*}
&(i)& T(\vw_{1}) \qp \vw_{2} = 0, \\
&(ii)&   \gwii{T(\vw_{1}) \, \vw_{2} \, \vw_{3} \, \vw_{4}}
    = \gwii{(\vw_{1} \qp \vw_{2}) \, \vw_{3} \, \vw_{4}} \\
&(iii)&
    \nabla_{\vw_{1}} \,\, T(\vw_{2}) = T \left(\nabla_{\vw_{1}} \vw_{2}
        \right) - \vw_{1} \qp \vw_{2}.
\end{eqnarray*}
Any vector field $\vw$ has the following decomposition
\begin{equation} \label{eqn:WTW}
\vw = T^{k} (\tau_{-}^{k}(\vw)) + \sum_{i=0}^{k-1} T^{i}
(\overline{\tau_{-}^{i}(\vw)})
\end{equation}
where $k$ is any positive integer (cf \cite[Equation (26)]{L2}). This decomposition
is very useful when applying topological recursion relations.
In particular, we will frequently use the decomposition
\[ \vw = \bvw + T(\tau_{-}(\vw)) \]
and call this the {\it standard decomposition} of $\vw$.
We also note that
$\{T^{k}(\ve_{i}) \mid i=1, \ldots, N, k \geq 0 \}$
gives a frame for the tangent bundle of the big phase space. An immediate consequence of
Lemma~\ref{lem:deridem} is the following
\begin{cor} \label{cor:deridem}
For any vector field $\vw$,
\begin{eqnarray*}
&(i)& \nabla_{_{T(\vw)}} \ve_{i} = - \vw \qp \ve_{i}, \\
&(ii)&    \nabla_{_{T^{2}(\vw)}} \ve_{i} = 0, \\
&(iii)& [T(\vw), \, \ve_{i}] = - T(\nabla_{\ve_{i}} \vw).
\end{eqnarray*}
for each $i$.
\end{cor}
{\bf Proof}: (i) and (ii) follow from Lemma~\ref{lem:deridem}
and above properties for $T$. (iii) follows from
(i) and \cite[Lemma 1.5]{L2}. Note that (iii) implies that
the frame $\{T^{k}(\ve_{i}) \mid i=1, \ldots, N, k \geq 0 \}$ is not commutative.
$\Box$

Motivated by
\[ \ve_{j} = \gwii{\ve_{j} \, \ve_{j} \, \gua} \ga, \]
we define
\[ \vf_{j} := \gwii{\ve_{j} \, \ve_{j} \, \ve_{j} \, \gua} \ga \]
for each $j = 1, \cdots, N$.
Then we have
\begin{lem} \label{lem:derejfj}
\begin{eqnarray*}
& (a) & \nabla_{\ve_{j}} \ve_{i} = \delta_{ij} \vf_{j} -
        \vf_{j} \qp \ve_{i} -   \vf_{i} \qp \ve_{j} \\
& (b) & \gwii{\ve_{i} \, \ve_{i} \, \ve_{j} \, \gua} \ga = \delta_{ij} \vf_{i}
   + \vf_{j} \qp \ve_{i} - \vf_{i} \qp \ve_{j} \\
& (c) & \gwii{\ve_{i} \, \ve_{j} \, \ve_{k} \, \vw}  = 0
        \hspace{20pt} {\rm for} \hspace{10pt} i \neq j  \neq k
        \hspace{20pt} {\rm and \, \, \,  any
        \,\,\, vector \, \, \, field}
        \hspace{10pt} \vw.
\end{eqnarray*}
\end{lem}
{\bf Proof}:
Since $\ve_{i} \qp \ve_{j} = 0$ for $i \ne j$, by \eqref{eqn:DerProd},
\[ 0 = \nabla_{_{\vw}} (\ve_{i} \qp \ve_{j})
    = (\nabla_{\vw} \ve_{i}) \qp \ve_{j}
        + (\nabla_{\vw} \ve_{j}) \qp \ve_{i} +
        \gwii{\vw \, \ve_{i} \, \ve_{j} \gua} \ga. \]
By Lemma~\ref{lem:deridem}, we have
\begin{equation} \label{eqn:weiej4pt}
\gwii{\vw \, \ve_{i} \, \ve_{j} \gua} \ga =
    -  \gwii{\vw \, \ve_{i} \, \ve_{i} \gua} \ga \qp \ve_{j}
        - \gwii{\vw \, \ve_{j} \, \ve_{j} \gua} \ga \qp \ve_{i}
\end{equation}
for any vector field $\vw$.
In particular, for $\vw = \ve_{i}$, we have
\begin{equation} \label{eqn:eieiej}
 \gwii{\ve_{i} \, \ve_{i} \, \ve_{j} \gua} \ga =
    -  \vf_{i} \qp \ve_{j}
        - \gwii{\ve_{i} \, \ve_{j} \, \ve_{j} \gua} \ga \qp \ve_{i}.
\end{equation}
Therefore
\[ \gwii{\ve_{i} \, \ve_{i} \, \ve_{j} \gua} \ga \qp \ve_{j} =
    -  \vf_{i} \qp \ve_{j}.
\]
Interchanging $i$ and $j$ in this formula then plugging in \eqref{eqn:eieiej},
we obtain (b). Lemma~\ref{lem:deridem} and (b) imply (a).
Replacing $\vw$ in \eqref{eqn:weiej4pt}
by $\ve_{k}$ for $k \neq i$ and $k \neq j$, then applying (b), we
obtain $\gwii{\ve_{i} \, \ve_{j} \, \ve_{k} \gua} \ga = 0$.
Since $\{ \ga \mid \alpha = 1, \ldots, N \}$ are linearly independent,
we have $\gwii{\ve_{i} \, \ve_{j} \, \ve_{k} \gua} = 0$.
Since  $\{ \gua \mid \alpha = 1, \ldots, N \}$ form a basis of the space
of primary vector fields, this proves (c) for the case where $\vw$ is a
primary vector field. On the other hand, if $\vw = T(\vv)$, then
(c) follows from the above mentioned properties for $T$ . Since any vector field can be
decomposed as a sum of a primary vector field and a vector field
of type $T(\vv)$, this proves (c).
$\Box$

Besides Lemma~\ref{lem:derejfj}, we also have the following
properties for the genus-0 4-point functions:
\begin{lem} \label{lem:4ptg0}
For $i \neq j$ and any vector fields $\vw$, $\vv$,
\begin{eqnarray*}
&(a)& \gwii{ \vw \, \vv \, \ve_{i} \, \ve_{j}} =
    - \gwii{ \vw \, (\vv \qp \ve_{j}) \, \ve_{i} \, \ve_{i}}
    - \gwii{ \vw \, (\vv \qp \ve_{i}) \, \ve_{j} \, \ve_{j}} \\
&(b)& \gwii{ \vw \, \ve_{i} \, \ve_{i} \, \ve_{j}}
        = - \gwii{ \vw \, \ve_{j} \, \ve_{j} \, \ve_{i}} \\
&(c)& \gwii{ \ve_{i} \, \ve_{i} \, \ve_{i} \, \ve_{j}}
    = \gwii{ \ve_{j} \, \ve_{j} \, \ve_{j} \, \ve_{i}}
    = - \gwii{ \ve_{i} \, \ve_{i} \, \ve_{j} \, \ve_{j}}, \\
&(d)& \gwii{\vs \, \vw \, \ve_{i} \, \ve_{j}} = 0.
\end{eqnarray*}
\end{lem}
{\bf Proof}: (a) is equivalent to \eqref{eqn:weiej4pt}. (b) is
obtained from (a) by setting $\vv = \ve_{i}$. (c) follows from
(b) by setting $\vw = \ve_{i}$ and $\vw = \ve_{j}$ respectively.
(d) follows from (a) with $\vv = \vs$ and (b). $\Box$

Now we look at some consequences of Lemma~\ref{lem:derejfj}.
\begin{cor} \label{cor:eiejcomm}
\[ [ \ve_{i}, \,  \ve_{j}] = 0. \]
\end{cor}
{\bf Proof}: This follows from Lemma~\ref{lem:derejfj} (a) and
the fact that
$[\ve_{i}, \, \ve_{j}] = \nabla_{\ve_{i}} \ve_{j} - \nabla_{\ve_{j}} \ve_{i}$.
$\Box$

\begin{cor} \label{cor:Fi}
\[ \vf_{j} = - \sum_{i=1}^{N}
    \left( \nabla_{\ve_{j}} \ve_{i} \right) \qp \ve_{i} \]
\end{cor}
{\bf Proof}: The two sides of this equation are equal
when multiplied by $\ve_{k}$ for all $k$ because of
 Lemma~\ref{lem:derejfj} (a).
$\Box$

Similar to \eqref{eqn:bvs}, we have
\begin{lem} \label{lem:tau-S}
\[\overline{\tau_{-}(\vs)} = - \sum_{i=1}^{N} \vf_{i}. \]
\end{lem}
{\bf Proof}:
By Lemma~\ref{lem:derejfj} (a),
\[ \nabla_{\bvs} \bvs = \sum_{i, j = 1}^{N} \nabla_{\ve_{i}} \ve_{j}
    = - \sum_{i=1}^{N} \vf_{i}. \]
On the other hand, by \cite[Lemma 1.8 (3)]{L2},
\[ \nabla_{\bvs} \bvs = \left( \bvs \gwii{\vs \, \vs \gua} \right) \ga
    = \bvs \qp \tau_{-}(\vs) = \overline{\tau_{-}(\vs)}. \]
The lemma follows.
$\Box$

For any vector field $\vw = \sum_{n, \alpha} f_{n, \alpha} \tau_{n}(\gamma_{\alpha})$,
define
\[ \vg * \vw = \sum_{n, \alpha} (n + b_{\alpha}) f_{n, \alpha} \tau_{n}(\gamma_{\alpha}). \]
This operator was used in \cite{L2} to give a recursive description for the Virasoro vector
fields.
\begin{lem} \label{lem:derXei}
\[ \nabla_{\vx} \ve_{i} = - \vg *  \ve_{i} + b_{1} \ve_{i},
    \hspace{20pt}
    \nabla_{\ve_{i}} \vx  = - \vg *  \ve_{i} + (b_{1}+1) \ve_{i}. \]
\end{lem}
{\bf Proof}:
By \cite[Equation (39)]{L2},
\begin{equation} \label{eqn:EulerEi4pt}
 \gwii{\vx \, \ve_{i} \, \ve_{i} \, \gua} \ga
    =  2 (\vg * \ve_{i}) \qp \ve_{i} - \vg * \ve_{i} - b_{1} \ve_{i}.
\end{equation}
The first formula follows from this equation and Lemma~\ref{lem:deridem}.
The second formula is a special case of a more general fact
(see the proof of \cite[Lemma 3.13]{L2}):
\[ \nabla_{\bvw} \vx = - \vg * \bvw + (b_{1}+1) \bvw \]
for any primary vector field $\bvw$.
$\Box$

\begin{cor} \label{cor:derXei}
\[ [ \vx, \, \ve_{i}] =[ \bvx, \, \ve_{i}] = - \ve_{i}. \]
\end{cor}
{\bf Proof}:
By Lemma~\ref{lem:derXei},
\[ [ \vx, \, \ve_{i}] = \nabla_{\vx} \ve_{i} - \nabla_{\ve_{i}} \vx
    = -\ve_{i}. \]
By Corollary~\ref{cor:deridem}, for any vector
field $\vw$,
\begin{equation} \label{eqn:bracketTWei}
[T(\vw), \, \ve_{i}]  =  \nabla_{T(\vw)} \ve_{i} - \nabla_{\ve_{i}} T(\vw)
= - T(\nabla_{\ve_{i}} \vw).
\end{equation}
Since $\tau_{-}$ commute with $\nabla$, by Lemma~\ref{lem:derXei},
\[ \nabla_{\ve_{i}} \, \tau_{-} (\vx) = \tau_{-} (\nabla_{\ve_{i}} \vx) = 0.\]
Therefore \eqref{eqn:bracketTWei} implies
\[ [T(\tau_{-}(\vx)), \, \ve_{i}] = 0. \]
The corollary then follows from the standard decomposition
 $\vx = \bvx + T(\tau_{-}(\vx))$.
$\Box$

\begin{cor}
\[ \ve_{j} \, u_{i} = \delta_{ij}. \]
Therefore in a vague sense,
we can think of $\ve_{i}$ as $\frac{\partial}{\partial u_{i}}$.
\end{cor}
{\bf Proof}: Since $\bvx = \sum_{j=1}^{N} u_{j} \ve_{j}$, by
Corollary~\ref{cor:eiejcomm},
\[ [\bvx, \, \ve_{i}] = - \sum_{j=1}^{N} (\ve_{i} u_{j}) \ve_{j}. \]
The desired formula then follows from Corollary~\ref{cor:derXei}.
$\Box$

\begin{lem} \label{lem:TWui}
 For any vector field $\vw$,
\[ T(\vw) \, u_{i} = 0. \]
\end{lem}
{\bf Proof}: By \eqref{eqn:bvx} and Corollary~\ref{cor:deridem},
\[ \nabla_{T(\vw)} \bvx = \sum_{i=1}^{N}\, \,  (T(\vw) \, u_{i}) \, \ve_{i}
        - u_{i} \, \vw \qp \ve_{i}. \]
On the other hand, \cite[Lemma 3.4]{L2} implies
\[  \nabla_{T(\vw)} \bvx = - \vw \qp \bvx =
    - \sum_{i=1}^{N}  u_{i} \, \vw \qp \ve_{i}.
\]
Subtracting these two equations, we obtain the desired formula.
$\Box$

Since for any vector field $\vw$, $\vw= \bvw + T(\tau_{-}(\vw))$,
we have
\begin{cor} \label{cor:Wui}
For all vector field $\vw$ on the big phase space,
\[ \vw u_{i} = \bvw u_{i}. \]
\end{cor}

\begin{lem} \label{lem:g*eij}
For $i \neq j$,
\[ (\vg * \ve_{i}) \qp \ve_{j} = (u_{j}-u_{i}) \vf_{i} \qp \ve_{j}. \]
\end{lem}
{\bf Proof}:
Multiplying both sides of \eqref{eqn:EulerEi4pt}
 by $\ve_{j}$ for $j \neq i$, we have
\[ (\vg * \ve_{i}) \qp \ve_{j}
    = - \gwii{\vx \, \ve_{i} \, \ve_{i} \, \gua} \ga \qp \ve_{j}. \]
On the other hand, using standard decomposition of $\vx$ and property of $T$, we have
\begin{eqnarray*}
 \gwii{\vx \, \ve_{i} \, \ve_{i} \, \gua} \ga
&=& \gwii{\bvx \, \ve_{i} \, \ve_{i} \, \gua} \ga
    + \gwii{T(\tau_{-}(\vx)) \, \ve_{i} \, \ve_{i} \, \gua} \ga \\
&=& \sum_{k=1}^{N} u_{k} \gwii{\ve_{k} \, \ve_{i} \, \ve_{i} \, \gua} \ga
    + \tau_{-}(\vx) \qp \ve_{i}.
\end{eqnarray*}
The Lemma then follows from Lemma~\ref{lem:derejfj}.
$\Box$

To completely determine $\vg * \ve_{i}$, it is convenient to introduce the
following bilinear form on the space of vector fields:
For any vector fields $\vw$ and $\vv$ on the big phase space, define
\[ <\vw, \vv> := \gwii{\vs \, \vw \, \vv} \]
The genus-0 topological recursion relation implies the following equation for genus-0 3-point
function
\[ \gwii{(\vw_{1} \qp \vw_{2}) \, \vw_{3} \, \vw_{4}}
    =\gwii{\vw_{1} \, (\vw_{2} \qp \vw_{3}) \, \vw_{4}}, \]
which implies the associativity of the quantum product. Therefore this bilinear form
has the following equivalent expressions
\begin{equation} \label{eqn:inner}
<\vw, \vv> = \gwii{\vs \, \vw \, \vv} = \gwii{\bvs \, \vw \, \vv}
    = \gwii{\vs \, \vs \, (\vw \qp \vv)}
    = (\vw \qp \vv) \gwii{\vs \, \vs}.
\end{equation}
The last equality follows from the fact that $\nabla_{\vw} \vs = 0$ if $\vw$ is a primary
vector field. In particular, this bilinear form is compatible with quantum product in the
following sense:
\[ <(\vw_{1} \qp \vw_{2}), \, \vw_{3}> = < \vw_{2}, \, (\vw_{1} \qp \vw_{3})>, \]
 and
\[ <\vw_{1}, \vw_{2}> = <\bvw_{1}, \bvw_{2}> \]
for any vector fields $\vw_{i}$.
The string equation implies that $<\ga, \gb> = \eta_{\alpha \beta}$ for
all $\alpha$ and $\beta$. So this bilinear form generalizes the Poincare
metric on the small phase space. But one should note that on the big phase
space, this bilinear form is highly degenerate because
\[ <T(\vw), \vv> = \gwii{\vs \,\, T(\vw) \, \vv} = 0 \]
for any vector fields $\vw$ and $\vv$. We also note that
by \cite[Lemma 1.8 and Equation (12)]{L2},
\begin{eqnarray}
 \vw_{1} <\vw_{2}, \vw_{3}>
&=& < \left\{ \nabla_{\vw_{1}} \vw_{2}
        + \vw_{1} \qp \tau_{-}(\vw_{2}) \right\}, \vw_{3}>  \nonumber \\
&&    + <\vw_{2},  \left\{ \nabla_{\vw_{1}} \vw_{3}
                    + \vw_{1} \qp \tau_{-}(\vw_{3}) \right\}>
                    \label{eqn:innermetric}
\end{eqnarray}
for any vector fields $\vw_{1}$, $\vw_{2}$, and $\vw_{3}$.
So this bilinear form is not compatible with the flat connection $\nabla$
in general. However, when restricted to the bundle of primary vector fields
over the big phase space, this bilinear form is non-degenerate and
compatible with $\nabla$, i.e.
\begin{equation} \label{eqn:CompConn}
 \vw_{1} <\vw_{2}, \vw_{3}>
= < \nabla_{\vw_{1}} \vw_{2}, \vw_{3}>
    + <\vw_{2},  \nabla_{\vw_{1}} \vw_{3}>
\end{equation}
if $\vw_{2}$ and $\vw_{3}$ are primary vector fields, and $\vw_{1}$
is an arbitrary vector field.
Note that by the last equality of \eqref{eqn:inner}, we have
\[ < \ve_{i}, \ve_{j}> = 0 \]
if $i \neq j$.
Therefore for any primary vector field $\vw$,
\begin{equation} \label{eqn:decomp-pvf}
\vw = \sum_{i=1}^{N} \frac{<\vw, \, \ve_{i}>}{<\ve_{i}, \ve_{i}>} \, \ve_{i}
\end{equation}
and
\begin{equation} \label{eqn:Wprodei}
\vw \qp \ve_{i}
=  \frac{<\vw, \, \ve_{i}>}{<\ve_{i}, \ve_{i}>} \, \ve_{i}
\end{equation}
for all $i$.

On the space of primary vector fields, the operator $\vg *$
has the following nice property:
\begin{equation} \label{eqn:G*inner}
 <\vg * \vw, \, \vv> + <\vw, \, \vg * \vv> = <\vw, \, \vv>
\end{equation}
for all primary vector fields $\vv$ and $\vw$. Since
$\{ \ga \mid \alpha=1, \ldots, N\}$ span the space of primary vector fields,
it suffices to check \eqref{eqn:G*inner} for $\vw= \ga$ and $\vv = \gb$
for arbitrary $\alpha$ and $\beta$. In this case we have
\[ <\vg * \ga, \, \gb> + <\ga, \, \vg * \gb>
    = (b_{\alpha} + b_{\beta}) \eta_{\alpha \beta}
    =  \eta_{\alpha \beta} = <\ga, \, \gb>.\]
An immediate consequence of \eqref{eqn:G*inner} is that
\begin{equation} \label{eqn:G*ei-inner}
< \vg * \ve_{i}, \, \ve_{i}> = \frac{1}{2} <\ve_{i}, \, \ve_{i}>
\end{equation}
for any $i$.

\begin{lem} \label{lem:g*ei}
\[ \vg * \ve_{i} = \frac{1}{2} \, \ve_{i}
        - u_{i} \vf_{i} + \vx \qp \vf_{i}. \]
\end{lem}
{\bf Proof}: By \eqref{eqn:G*ei-inner} and
\eqref{eqn:Wprodei},
\[ (\vg * \ve_{i}) \qp \ve_{i} = \frac{1}{2} \, \ve_{i}. \]
 For $j \neq i$,
$(\vg * \ve_{i}) \qp \ve_{j}$ can be computed
using Lemma~\ref{lem:g*eij}.
The two sides  of the equation in this lemma are equal when multiplied by
$\ve_{j}$ for all $j$. Since both sides are primary vector fields, this
proves the lemma.
$\Box$

\section{Rotation coefficients on the big phase space}
\label{sec:rotcoeff}

Rotation coefficients played a very important role in the study of semisimple Frobenius manifolds
by Dubrovin (cf. \cite{D}). In this section, we define rotation coefficients for
the quantum product on the big phase space and study their basic properties.
Genus-0 correlation functions obey the WDVV equation and its derivatives.
Derivatives of WDVV equations are not symmetric in general. This may cause annoying troubles
in studying relations among complicated expressions of genus-0 invariants (See, for example,
\cite{L3} for the occurrence of derivatives of WDVV equation among certain genus-1
equations). However, in the semisimple case, such genus-0 relations are somehow trivialized
after introducing rotation coefficients.
Many formulas in this section which only involve primary fields are analogue
(in a slightly more complicated way) of the corresponding formulas
on the small phase space as given in \cite{D}. We will also discuss the behavior of relevant
functions under derivatives along descendant vector fields.

For any primary vector field $\vw$,
define
\[ \| \vw \| : = \sqrt{<\vw, \vw>}. \]
For each $i$, define
\[ g_{i} := \| \ve_{i} \|^{2} = <\ve_{i}, \ve_{i}>. \]
Since $ < \ve_{i}, \ve_{j}> = 0 $ if $i \neq j$,
functions $g_{1}, \ldots, g_{N}$ completely
determines $<\cdot, \cdot>$ in the semisimple case.

\begin{lem} \label{lem:dergij}
For all $i$ and any vector field $\vw$,
\begin{eqnarray*}
&(a)& \ve_{j} g_{i} = \gwii{\vs \, \vs \,
    \left\{ \tau_{-}(\vs) \qp \ve_{i} \qp \ve_{j} +
            \nabla_{\ve_{j}} \ve_{i} \right\}}
            \,\,\, {\rm for \,\,\, all} \,\,\, j,
     \\
&(b)& T(\vw) \, g_{i} = - 2 <\vw, \ve_{i}> \\
&(c)& \vs g_{i} = 0, \\
&(d)& \vx g_{i} = (2 b_{1} - 1) g_{i}.
\end{eqnarray*}
\end{lem}
{\bf Proof}:
Since $\nabla_{\ve_{j}} \vs = - \tau_{-}(\ve_{j}) = 0$ and
$g_{j} = \gwii{\vs \, \vs \, \ve_{i}}$,
\[ \ve_{j} g_{i} = \gwii{\vs \, \vs \, \ve_{i} \, \ve_{j}}
        + \gwii{\vs \, \vs \, (\nabla_{\ve_{j}} \ve_{i})}. \]
By the string equation, this implies (a). Since
$\nabla_{T(\vw)} \vs = - \tau_{-}T(\vw) = - \vw$,
(b) follows from \cite[Corollary 1.6]{L2} and Corollary~\ref{cor:deridem}.
By Lemma~\ref{lem:deridem} and the string equation,
\[ \nabla_{\vs} \, \, \ve_{i}
    = R_{i}(\gwii{\vs \, \ve_{i} \, \ve_{i} \, \gua} \ga) = 0. \]
So (c) follows from \eqref{eqn:CompConn}
since $g_{i} = <\ve_{i}, \ve_{i}>$.
(d) follows from Lemma~\ref{lem:derXei} and \eqref{eqn:G*ei-inner}.
$\Box$

Define {\it Rotation coefficients} by
\[ r_{ij} := \frac{\ve_{j}}{\sqrt{g_{j}}} \sqrt{g_{i}}
    = \frac{ \ve_{j} \, g_{i}}{2 \, \sqrt{g_{i} g_{j}}}. \]
Covariant derivatives of idempotents can be computed using rotation
coefficients.
\begin{lem} \label{lem:rotder}
\begin{eqnarray*}
 \nabla_{\ve_{i}} \ve_{j}
&=& r_{ij} \left( \sqrt{\frac{g_{j}}{g_{i}}} \, \ve_{i}
            + \sqrt{\frac{g_{i}}{g_{j}}} \, \ve_{j} \right)
  - \delta_{ij} \sum_{k=1}^{N}
         r_{ik} \,  \sqrt{\frac{g_{i}}{g_{k}}} \,\, \ve_{k}
\end{eqnarray*}
\end{lem}
{\bf Proof}:
If $i \neq j$,
\[ <\ve_{i}, \, \nabla_{\ve_{j}} \ve_{i}> =
    \frac{1}{2} \ve_{j} \, g_{i} = r_{ij} \, \sqrt{g_{i} g_{j}} \]
\[ <\ve_{j}, \, \nabla_{\ve_{j}} \ve_{i}> =
    <\ve_{j}, \, \nabla_{\ve_{i}} \ve_{j}> =
    \frac{1}{2} \ve_{i} \,  g_{j} = r_{ij} \, \sqrt{g_{i} g_{j}}. \]
The formula then follows from the fact that $<\ve_{i}, \ve_{j}> = g_{i} \delta_{ij}$.

For $k \neq i$, since $<\ve_{k}, \, \ve_{i}> = 0$,
\[ <\ve_{k}, \, \nabla_{\ve_{i}} \ve_{i}>
    = - <\nabla_{\ve_{i}} \ve_{k}, \,  \ve_{i}>
    = -  r_{ik} \, \sqrt{g_{i} g_{k}}.\]
Moreover
\[ <\ve_{i}, \, \nabla_{\ve_{i}} \ve_{i}> =
    \frac{1}{2} \ve_{i} \, g_{i} = r_{ii} \, g_{i}.  \]
So
\[ \nabla_{\ve_{i}} \ve_{i}
=  r_{ii} \, \ve_{i} - \sum_{k \neq i}
         r_{ik} \,  \sqrt{\frac{g_{i}}{g_{k}}} \,\, \ve_{k}. \]
This proves the lemma.
$\Box$

By Corollary~\ref{cor:Fi}, we have
\begin{cor} \label{cor:FiEj}
For every $i$,
\[ \vf_{i} = - \sum_{j=1}^{N}
         r_{ij} \,  \sqrt{\frac{g_{i}}{g_{j}}} \,\, \ve_{j}. \]
\end{cor}

Together with Lemma~\ref{lem:g*ei}, this corollary implies
\begin{cor} \label{cor:G*rot}
\[ \vg * \ve_{i} = \frac{1}{2} \ve_{i}
    + \sum_{j} (u_{i} - u_{j}) r_{ij} \, \sqrt{\frac{g_{i}}{g_{j}}}
    \, \ve_{j}
\]
for all $i$.
\end{cor}

More properties of rotation coefficients are collected in the following:
\begin{lem} \label{lem:derRot}
For $i, j = 1, \ldots, N$ and any vector field $\vw$,
\begin{eqnarray*}
&(a)& r_{ij} = r_{j \, i}, \\
&(b)& T(\vw) \, r_{ij} = \delta_{ij} \left\{  - \frac{<\tau_{-}(\vw), \, \ve_{i}>}{g_{i}}
        + \sum_{k=1}^{N} \frac{r_{ik}}{\sqrt{g_{i} g_{k}}}
        \left<\vw, \,\,\,  \ve_{k} \right>
         \right\}, \\
&(c)& \vs \, r_{ij} = 0, \\
&(d)& \vx \, r_{ij} = - r_{ij}.
\end{eqnarray*}
\end{lem}
{\bf Proof}:
(a) follows from Lemma~\ref{lem:dergij} (a) since
$\nabla_{\ve_{i}} \ve_{j} = \nabla_{\ve_{j}} \ve_{i}$.

To prove (b), we first note that
\[ T(\vw) r_{ij}
    = T(\vw) \left( \frac{\ve_{j} g_{i}}{2 \sqrt{ g_{i}g_{j}}} \right)
    = \frac{ T(\vw) \ve_{j} g_{i}}{2 \sqrt{ g_{i}g_{j}}}
    + (\ve_{j} g_{i}) T(\vw)
        \left( \frac{1}{2 \sqrt{ g_{i}g_{j}}} \right). \]
By Lemma~\ref{lem:dergij} (b), the second term is
\[ (\ve_{j} g_{i}) T(\vw)
        \left( \frac{1}{2 \sqrt{ g_{i}g_{j}}} \right) =
     r_{ij} \left( \frac{<\vw, \ve_{i}>}{g_{i}} +
                    \frac{<\vw, \ve_{j}>}{g_{j}} \right). \]
By Corollary~\ref{cor:deridem} (iii) and Lemma~\ref{lem:dergij} (b),
the first term is
\begin{eqnarray*}
\frac{ T(\vw) \ve_{j} g_{i}}{2 \sqrt{ g_{i}g_{j}}}
&=& \frac{1}{2 \sqrt{ g_{i}g_{j}}}
\left( \ve_{j} T(\vw) g_{i} - T(\nabla_{\ve_{j}} \vw) g_{i} \right) \\
&=& \frac{1}{2 \sqrt{ g_{i}g_{j}}}
\left( -2 \ve_{j} <\vw, \, \ve_{i}> + 2 <\nabla_{\ve_{j}} \vw, \, \ve_{i}>
    \right).
\end{eqnarray*}
By \eqref{eqn:innermetric},
\begin{eqnarray*}
 \ve_{j} <\vw, \, \ve_{i}>
    &=& <\nabla_{\ve_{j}} \vw + \ve_{j} \qp \tau_{-}(\vw), \, \ve_{i}>
        + <\vw, \nabla_{\ve_{j}} \ve_{i}>  \\
    &=& <\nabla_{\ve_{j}} \vw, \, \ve_{i}>
        + <\vw, \nabla_{\ve_{j}} \ve_{i}>
        + \delta_{ij}<\tau_{-}(\vw), \ve_{i}>.
\end{eqnarray*}
Therefore we have
\[ \frac{ T(\vw) \ve_{j} g_{i}}{2 \sqrt{ g_{i}g_{j}}}
    = - \, \frac{<\vw, \, \nabla_{\ve_{j}} \ve_{i}>
        + \delta_{ij} <\tau_{-}(\vw), \ve_{i}>}{ \sqrt{ g_{i}g_{j}}}.
        \]
(b) then follows from Lemma~\ref{lem:rotder}.

The string equation implies that
$\gwii{\vs \, \ve_{i} \, \ve_{i} \, \gua} = 0$ for any $i$ and $\alpha$
(cf. \cite[Equation (12)]{L2}).
By Lemma~\ref{lem:deridem},
\[ \nabla_{\vs} \ve_{i} =  0 \]
for any $i$.
Since $\nabla_{\ve_{i}} \vs = - \tau_{-}(\ve_{i}) = 0$, we have
$[\vs, \, \ve_{i}] = 0$ for all $i$. Therefore (c) follows
from Lemma~\ref{lem:dergij} (c).

By Lemma~\ref{lem:derXei}, $[\vx, \, \ve_{j}] = - \ve_{j}$, so
\[ \vx r_{ij} = \vx \left( \frac{\ve_{j} g_{i}}{2 \sqrt{ g_{i}g_{j}}} \right)
    = \frac{\ve_{j} \vx g_{i} - \ve_{j} g_{i}}{2 \sqrt{ g_{i}g_{j}}}
    + (\ve_{j} g_{i}) \vx
        \left( \frac{1}{2 \sqrt{ g_{i}g_{j}}} \right). \]
Applying Lemma~\ref{lem:dergij} (d) and replacing $\ve_{j} g_{i}$
by $2 r_{ij} \sqrt{ g_{i}g_{j}}$, we obtain (d).
$\Box$

{\bf Remark}: For the dilaton vector field $\vd = T(\vs)$,
\[ \vd \, r_{ij} = 0\]
for all $i$ and $j$. So for $\vl_{o} = -\vx - (b_{1}+1) \vd$ and
any $i$, $j$,
\[ \vl_{0} r_{ij} = r_{ij}. \]

Derivatives of rotation coefficients along primary vector fields can be
computed using the following formulas:
\begin{lem} \label{lem:derRotPrim}
For $i \neq j$,
\begin{eqnarray*}
&(a)& \ve_{k} \, r_{ij} = r_{ik} \, r_{jk} \,\,\,
        {\rm if} \,\,\, i, j, k \,\,\, {\rm are \,\,\, distinct}, \\
&(b)&  \ve_{i} \, r_{ij} = \frac{1}{u_{j}-u_{i}}
    \left\{ r_{ij} + \sum_{k \neq i, j} (u_{k}-u_{j}) r_{ik} r_{jk}
    \right\}, \\
&(c)& \ve_{j} \, r_{ii} = r_{ij}^{2} +
    \sqrt{\frac{g_{j}}{g_{i}}} (\ve_{i} \, r_{ij}
            - r_{ii} r_{ij}), \\
&(d)& \ve_{i} \, r_{ii} = - r_{ii}^{2} -  \sum_{j \neq i}
    \left\{ 2 r_{ij}^{2} +  \sqrt{\frac{g_{j}}{g_{i}}} \, \
        r_{ij} r_{jj} - \sqrt{\frac{g_{j}}{g_{i}}} \, \
        \ve_{j} r_{ij} \right\}
         + \frac{1}{g_{i}} <\tau_{-}^{2}(\vs), \, \ve_{i}>.
\end{eqnarray*}
\end{lem}
{\bf Proof}:
By definition of rotation coefficients,
$\ve_{i} \, \sqrt{g_{j}} = \sqrt{g_{i}} \, r_{ij}$ for all $i$ and $j$.
So
\[ \ve_{i} \sqrt{\frac{g_{j}}{g_{k}}}
    =  \sqrt{\frac{g_{i}}{g_{k}}} \, r_{ij} -
        \frac{\sqrt{g_{i} g_{j}}}{g_{k}} \, r_{ik} \]
for all $i$, $j$, and $k$.
Use this formula and Lemma~\ref{lem:rotder}, we obtain
\[ \left< \nabla_{\ve_{i}} \nabla_{\ve_{j}} \ve_{k}, \, \ve_{j} \right>
    = r_{jk} r_{ik} \sqrt{g_{i}g_{j}}
        + \sqrt{g_{j} g_{k}} \, \, \ve_{i} r_{jk} \]
and
\[\left< \nabla_{\ve_{j}} \nabla_{\ve_{i}} \ve_{k}, \, \ve_{j} \right>
  = r_{ik} \left\{ r_{jk} \sqrt{g_{i}g_{j}} +
        r_{ij} \sqrt{g_{j}g_{k}} \right\} \]
for distinct $i$, $j$, and $k$. Since $\nabla$ is a flat connection,
the left hand side of these two equations are equal. The equality of the
right hand sides of these two equations is precisely the formula in (a).

By Lemma~\ref{lem:derRot} (b) and (c),
\begin{equation} \label{eqn:bSrij}
 \bvs \, r_{ij} = 0 \hspace{20pt} {\rm if} \, \, \, i \neq j,
\end{equation}
Since $\bvs=\sum_{k=1}^{N} \ve_{k}$, the formula in (a)  implies
\begin{equation} \label{eqn:Eirij}
 \ve_{i} \, r_{ij} + \ve_{j} r_{ij} = - \sum_{k \neq i, j} r_{ik} r_{jk}
\end{equation}
for $i \neq j$.
Moreover, by Lemma~\ref{lem:derRot} (b) and (d),
\begin{equation} \label{eqn:bXrij}
 \bvx \, r_{ij} = - r_{ij} \hspace{20pt} {\rm if} \, \, \, i \neq j.
\end{equation}
Since $\bvx = \sum_{i=1}^{N} u_{i} \ve_{i}$, using (a) again, we obtain
\begin{equation} \label{eqn:uEirij}
u_{i} \ve_{i} r_{ij} + u_{j} \ve_{j} r_{ij} = -r_{ij}
    - \sum_{k \neq i, j} u_{k} r_{ik} r_{jk}
\end{equation}
for $i \neq j$.
Solving $\ve_{i} r_{ij}$ from \eqref{eqn:Eirij} and \eqref{eqn:uEirij},
we obtain (b).

The formula in (c) follows directly from the definition of $r_{ij}$
and the fact  $[\ve_{i}, \, \ve_{j}] = 0$. We can also obtain (c)
and \eqref{eqn:Eirij} from the fact that $\nabla$ is a flat connection.
Together with (a), these are perhaps all what all we can get from the flatness
of $\nabla$.

By Lemma~\ref{lem:derRot} (b) and (c),
\[ \bvs \, r_{ii} = - T(\tau_{-}(\vs)) \, r_{ii}
    = - \left< \overline{\tau_{-}(\vs)}, \, \,
        \sum_{k=1}^{N} \frac{r_{ik}}{\sqrt{g_{i}g_{k}}} \ve_{k}
            \right> + \frac{1}{g_{i}} <\tau_{-}^{2}(\vs), \, \ve_{i}>.
\]
Then applying Lemma~\ref{lem:tau-S} and Corollary~\ref{cor:FiEj}, we obtain
\begin{equation} \label{eqn:bvSrii}
 \bvs \, r_{ii} = - \sum_{j, k=1}^{N} \sqrt{\frac{g_{j}}{g_{i}}}
        \, \, r_{ik} \, r_{jk}
        + \frac{1}{g_{i}} <\tau_{-}^{2}(\vs), \, \ve_{i}>.
\end{equation}
Since $\bvs = \sum_{j=1}^{N} \ve_{j}$, we can solve
$\ve_{i} \, r_{ii}$ from this equation. Note that
if $i \neq j$, (a) implies that
\[ \sum_{k=1}^{N} r_{ik} \, r_{jk} = (r_{ii}+ r_{jj})r_{ij}
        + \bvs \, r_{ij} - (\ve_{i} + \ve_{j}) r_{ij}.\]
We can simplify the formula for $\ve_{i} \, r_{ii}$ using this
equation and \eqref{eqn:bSrij}, then use (c) to obtain (d).
$\Box$

{\bf Remark}: This lemma allows us to compute all derivatives of
rotation coefficients except $\ve_{i} r_{ii}$ in terms of functions
$g_{i}$, $r_{ij}$, and $u_{i}$. For $\ve_{j} r_{ii}$ we have:
\begin{eqnarray}
\ve_{j} \, r_{ii} & = & r_{ij}^{2} +  \sqrt{\frac{g_{j}}{g_{i}}} \, \,
        \frac{1}{u_{j}-u_{i}}
    \left\{ r_{ij} + \sum_{k} (u_{k}-u_{j}) r_{ik} r_{jk}
    \right\} \, \, \, {\rm if} \, \, \, i \neq j.
\end{eqnarray}
On the small phase space, $<\tau_{-}^{2}(\vs), \, \ve_{i}> = 0$.
So we have
\[ \ve_{i} \, r_{ii}  =  - r_{ii}^{2} - 2 \sum_{j \neq i} r_{ij}^{2}
    - \sum_{j \neq i}
     \sqrt{\frac{g_{j}}{g_{i}}} \, \, \frac{1}{u_{j}-u_{i}}
    \left\{ r_{ij} + \sum_{k} (u_{k}-u_{i}) r_{ik} r_{jk}
    \right\} \]
on the small phase space. We also note that on the small phase space,
$\bvs = \vs$ and $\bvx = \vx$. So on the small phase space,
we have for every $i$,
\[ \sum_{j} \sqrt{g_{j}} \, r_{ij} = \bvs \sqrt{g_{i}} = 0 \]
and
\[ \sum_{j} u_{j} \sqrt{g_{j}} \, r_{ij} = \bvx \sqrt{g_{i}}
    = (b_{1} - \frac{1}{2}) \sqrt{g_{i}}.
\]
Therefore many formulas can be simplified on the small phase space.

On the big phase space $<\tau_{-}^{2}(\vs), \, \ve_{i}> \neq 0$ in general.
This term has to be included when we compute $\ve_{i} \, r_{ii}$.  This is a
typical big phase space phenomenon.   When we taking higher
order derivatives of rotation coefficients, we will also encounter
$<\tau_{-}^{k}(\vs), \, \ve_{i}>$ for $k \geq 2$. However for most purpose, these
terms do not affect calculation. This is true for the proof of genus-1 and genus-2
Virasoro conjecture for the semisimple case.

Lemma~\ref{lem:derRotPrim} (a) and (b) can be encoded in a nice
matrix equation as observed in \cite{D} for the small phase space.
Let $\Gamma = (r_{ij})_{N \times N}$, $U$ the diagonal
$N \times N$ matrix whose diagonal entries are $\{u_{1}, \ldots, u_{N}\}$,
and $E_{i}$ the $N \times N$ matrix with all entries equal to $0$ except
that the $i$-th entry along diagonal is equal to $1$. Define
\[ V := [\Gamma, U].\]
 Note that the $(i, j)$-th entry of $V$ is
\[ v_{ij} = (u_{j} -u_{i}) r_{ij}. \]
These functions appear naturally in the expression of $G*\ve_{i}$.
Corollary~\ref{cor:G*rot} can be rewritten as
\[ G*\ve_{i} = \frac{1}{2} \ve_{i}
    - \sum_{j} v_{ij} \sqrt{\frac{g_{i}}{g_{j}}} \ve_{j}. \]
Since $V$ is a skew symmetric matrix,
eigenvalues of $V$ appear in opposite pairs.
Lemma~\ref{lem:derRotPrim} (a) and (b)
are equivalent to
\begin{equation} \label{eqn:matrixDerRot}
\ve_{k} \, V = [V, \, [E_{k}, \Gamma]]
\end{equation}
for all $k$. It turns out $b_{\alpha} - \frac{1}{2}$ is an eigenvalue of
$V$ for any $\alpha$. The corresponding eigenvector is the column vector
$\psi_{\alpha} := (\psi_{1 \alpha}, \ldots, \psi_{N \alpha})^{T}$ where
\begin{equation} \label{eqn:psiia}
 \psi_{i \alpha} := \eta_{\alpha \beta} \frac{\ve_{i}}{\sqrt{g_{i}}}
            t_{0}^{\beta}
        =  \frac{\ve_{i}}{\sqrt{g_{i}}} \gwii{\vs \, \ga}
        = \frac{1}{\sqrt{g_{i}}} \gwii{\vs \, \ve_{i} \,  \ga}.
\end{equation}
The last two equalities are due to the string equation.
(Note that in the definition of $\psi_{\alpha}$,
if we replace $\eta^{\alpha \beta} t_{0}^{\beta}$ by $t_{0}^{\alpha}$
we still obtain an eigenvector of $V$ but with eigenvalue
$\frac{1}{2} - b_{\alpha}$.)
The last equality of \eqref{eqn:psiia}
also implies that
\begin{equation} \label{eqn:guatoei}
 \frac{\ve_{i}}{\sqrt{g_{i}}} = \psi_{i \alpha} \gua.
\end{equation}
Therefore Corollary~\ref{cor:G*rot} is equivalent to the fact that
\[ V \psi_{\alpha} = (b_{\alpha} - \frac{1}{2} )\psi_{\alpha} \]
for all $\alpha$. The $N \times N$ matrix
\[ \psi := (\psi_{i \alpha})_{N \times N} = (\psi_{1}, \ldots, \psi_{N}) \]
is the transition matrix from one frame of primary vector fields
$\{ \gua \mid \alpha=1, \ldots, N\}$ to another frame
$\{ \frac{\ve_{i}}{\sqrt{g_{i}}} \mid i = 1, \ldots, N \}$.
Since $\{ \frac{\ve_{i}}{\sqrt{g_{i}}} \mid i = 1, \ldots, N \}$ is an
orthonormal frame, we have
\begin{equation} \label{eqn:psi-inverse}
 \delta_{ij} = \left<\frac{\ve_{i}}{\sqrt{g_{i}}}, \,
        \frac{\ve_{j}}{\sqrt{g_{j}}} \right>
    = \psi_{i \alpha} \eta^{\alpha \beta} \psi_{j \beta}
\end{equation}
for any $i$ and $j$. This is equivalent to
\begin{equation} \label{eqn:psipsi-1}
 {\rm Id} = \psi \eta^{-1} \psi^{T}, \hspace{30pt}
 \psi^{-1} = \eta^{-1} \psi^{T}, \hspace{30pt} \psi^{T} \psi = \eta.
\end{equation}
Consequently, by \eqref{eqn:guatoei},
\begin{equation} \label{eqn:eitoga}
\ga = \sum_{i=1}^{N} \psi_{i \alpha} \frac{\ve_{i}}{\sqrt{g_{i}}}
\end{equation}
for all $\alpha$.

Note that quantum products of primary vector fields
are encoded in the transition
matrix $\psi$. In fact, by associativity,
\begin{equation} \label{eqn:3ptRot}
 \gwii{\ve_{i} \, \ve_{j} \, \ve_{k}} =
    \gwii{\vs \, \vs \, (\ve_{i} \qp \ve_{j} \qp \ve_{k})} =
    \delta_{ij} \delta_{ik} g_{i}
\end{equation}
for any $i, j, k$. By \eqref{eqn:eitoga},
\begin{equation} \label{eqn:3ptpsi}
 \gwii{\ga \, \gb \, \gm} = \sum_{i=1}^{N}
    \frac{ \psi_{i \alpha} \psi_{i \beta} \psi_{i \mu}}{\sqrt{g_{i}}}
\end{equation}
for all $\alpha$, $\beta$, and $\mu$.
On the small phase space, $\vs = \gamma_{1}$. So on the small phase
$\psi_{i 1} = \gwii{\vs \, \ve_{i} \, \gamma_{1}} / \sqrt{g_{i}}
= \sqrt{g_{i}}$. So when restricted to the small phase space,
\eqref{eqn:3ptpsi} coincides with the corresponding formula in
\cite{D}. But this is not true on the big phase space.

We can also express $V$ in terms of $\psi_{i \alpha}$. Let
$A$ be the diagonal matrix whose $\alpha$-th entry along the diagonal
is $b_{\alpha} - \frac{1}{2}$. Then $V \psi = \psi A$ since column vectors
of $\psi$ are eigenvectors of $V$ with eigenvalues equal to
diagonal entries of
$A$. Therefore we have
\[ V = \psi A \psi^{-1} = \psi A \eta^{-1} \psi^{T}\]
or equivalently
\begin{equation} \label{eqn:rotpsi}
 (u_{j} - u_{i}) r_{ij} = v_{ij}
    = \sum_{\alpha, \beta} (b_{\alpha} - \frac{1}{2})
        \psi_{i \alpha} \eta^{\alpha \beta} \psi_{j \beta}
\end{equation}
for all $i$ and $j$. This formula tells us how to compute rotation coefficients
in terms of $\psi$.

Derivatives of $\psi_{i \alpha}$ can be computed in the following way:
By \eqref{eqn:guatoei} and Lemma~\ref{lem:rotder},
\[ (\ve_{k} \psi_{i \alpha}) \gua
    = \nabla_{\ve_{k}} (\psi_{i \alpha} \gua) =
    \nabla_{\ve_{k}} \left(\frac{\ve_{i}}{\sqrt{g_{i}}} \right)
    = r_{ki} \frac{\ve_{k}}{\sqrt{g_{k}}}
    = r_{ki} \psi_{k \alpha} \gua \]
for $k \neq i$. So
\[ \ve_{k} \psi_{i \alpha} = r_{ki} \psi_{k \alpha} \]
for all $\alpha$ and $i \neq k$. Similar proof also shows that
\[ \ve_{i} \psi_{i \alpha} = - \sum_{j \neq i} r_{ij} \psi_{j \alpha} \]
for all $\alpha$ and $i$. These two formulas can be combined as an equation
for the vector $\psi_{\alpha}$:
\[ \ve_{k} \psi_{\alpha} = - [E_{k}, \, \Gamma] \psi_{\alpha} \]
for all $k$ and $\alpha$. Using \eqref{eqn:eitoga}, we also have
\[ \gb \, \psi_{i \alpha} = \sum_{j = 1}^{N} r_{ij} \psi_{j \alpha}
            \left( \frac{\psi_{j \beta}}{\sqrt{g_{j}}} -
                \frac{\psi_{i \beta}}{\sqrt{g_{i}}} \right) \]
for all $i$, $\alpha$ and $\beta$.

We can represent genus-0 $k$-point functions in terms of rotation coefficients
$r_{ij}$ and functions $g_{i}$ and $u_{i}$. To do this, we need generalize
\eqref{eqn:3ptRot} to $k$-point functions. This can be done by repeatedly
taking derivatives of \eqref{eqn:3ptRot} along idempotent vector fields
and applying Lemma~\ref{lem:rotder} and \ref{lem:derRotPrim}.
The formula becomes more and more complicated as $k$ becomes larger and larger.
For 4-point functions, we have
\begin{lem} \label{lem:4ptRot}
For distinct $i$, $j$, $k$, and any vector field $\vw$,
\begin{eqnarray*}
&(i)& \gwii{\ve_{i} \, \ve_{i} \, \ve_{i} \, \ve_{i}}
    = - g_{i} r_{ii}, \\
&(ii)& \gwii{\ve_{j} \, \ve_{i} \, \ve_{i} \, \ve_{i}}
    = - \gwii{\ve_{j} \, \ve_{j} \, \ve_{i} \, \ve_{i}}
    = - \sqrt{g_{i} g_{j}} \,  r_{ij}, \\
&(iii)& \gwii{\ve_{i} \, \ve_{j} \, \ve_{k} \, \vw} = 0.
\end{eqnarray*}
\end{lem}

\section{Applications to higher genus Gromov-Witten invariants }

Idempotents can be applied to study higher genus Gromov-Witten invariants.
For example, they can be used to solve universal equations and prove the
Virasoro conjecture for manifolds with semisimple quantum cohomology
up to genus-2. There is no doubt that they can also be applied to Gromov-Witten
invariants of genus bigger than 2 once the corresponding universal equations
are found. In this section, we will mainly illustrate how to apply idempotents to study
genus-1 Gromov-Witten invariants. We will briefly comment on the genus-2 case, but details
will be given in a separate paper \cite{L4}.

For any vector fields $v_{1}, \ldots v_{4}$ on the big
phase space, we define
\begin{eqnarray*}
G_{0}(v_{1}, v_{2}, v_{3}, v_{4}) & = &
    \sum_{g \in S_{4}} \sum_{\alpha, \beta} \left\{
        \frac{1}{6} \left<\left< v_{g(1)} v_{g(2)}v_{g(3)}
                    \gamma^{\alpha} \right>\right>_{0}
        \left<\left< \gamma_{\alpha} v_{g(4)} \gamma_{\beta}
                    \gamma^{\beta} \right>\right>_{0}
        \right. \\
    && \hspace{40pt}
    + \frac{1}{24} \left<\left< v_{g(1)} v_{g(2)}v_{g(3)} v_{g(4)}
                    \gamma^{\alpha} \right>\right>_{0}
    \left<\left< \gamma_{\alpha}  \gamma_{\beta}
                    \gamma^{\beta} \right>\right>_{0}
            \\
    && \hspace{40pt} \left.
    - \frac{1}{4} \left<\left< v_{g(1)} v_{g(2)}
            \gamma^{\alpha} \gamma^{\beta} \right>\right>_{0}
        \left<\left< \gamma_{\alpha}  \gamma_{\beta}
                v_{g(3)} v_{g(4)} \right>\right>_{0}
            \right\},
\end{eqnarray*}
and
\begin{eqnarray*}
G_{1}(v_{1}, v_{2}, v_{3}, v_{4}) & = &
    \sum_{g \in S_{4}}
        3 \left<\left< \{v_{g(1)} \qp v_{g(2)} \}
            \{v_{g(3)} \qp v_{g(4)} \}
                     \right>\right>_{1}
        \\
    &&
    - \sum_{g \in S_{4}}
    4 \left<\left< \{v_{g(1)} \qp v_{g(2)} \qp
            v_{g(3)} \} v_{g(4)}
                     \right>\right>_{1}
    \\
    &&
     - \sum_{g \in S_{4}} \sum_{\alpha}
        \left<\left< \{ v_{g(1)} \qp v_{g(2)} \}
                v_{g(3)} v_{g(4)}
                    \gamma^{\alpha} \right>\right>_{0}
    \left<\left< \gamma_{\alpha} \right>\right>_{1}
        \\
    &&
     + \sum_{g \in S_{4}} \sum_{\alpha}
         2 \left<\left< v_{g(1)}  v_{g(2)} v_{g(3)}
                    \gamma^{\alpha} \right>\right>_{0}
    \left<\left< \{\gamma_{\alpha} \qp v_{g(4)} \}
            \right>\right>_{1}
\end{eqnarray*}
where $S_{4}$ is the permutation group of 4 elements. Note that
$G_{0}$ is determined solely by genus-0 data, while each term in
$G_{1}$ contains genus-1 information. These two tensors are
connected by the following equation due to Getzler (cf. \cite{G1}):
\begin{equation} \label{eqn:G}
    G_{0} + G_{1} = 0.
\end{equation}
On the small phase space, when the quantum cohomology is semisimple,
Dubrovin and Zhang \cite{DZ1} solved genus-1 generating functions from this equation
up to a constant.
The generating function on the big phase space can be obtained by a so called
constitutive relation, which expresses the big phase space generating function
in terms of the small phase space generating function (with some extra terms).
We will deal with this equation directly on the big phase space, therefore combine
two steps into one step.
 We will also need to use
 the {\it genus-1 topological recursion relation} which has the following
form:
\begin{equation} \label{eqn:TRRg1}
\gwiione{T({\cal W})} \, \, = \, \,
 \frac{1}{24} \gwii{{\cal W} \, \gum  \, \gm}
\end{equation}
for any vector field $\vw$. Its first derivative has the form:
\begin{equation}
\gwiione{T({\cal W}) \, \vv} \, \, = \, \,
    \gwiione{ \{\vw \bullet \vv\} } +
     \frac{1}{24} \gwii{{\cal W} \, \vv \, \gum  \, \gm}
\label{eqn:derTRRg1}
\end{equation}
for all vector field $\vw$ and $\vv$ (cf. \cite[Equation
(14)]{L2}).

\begin{pro} \label{pro:g1ei}
For each $i$,
\[ \phi_{i}
    := \gwiione{\ve_{i}}
    = \frac{1}{24} \left( \gwii{\ve_{i} \, \tau_{-}(\vl_{0}) \, \ga \, \gua}
        - G_{0}(\ve_{i}, \ve_{i}, \ve_{i}, \bvx) \right),\]
where $\vl_{o} = -\vx - (b_{1}+1) T(\vs)$.
Note the right hand side of this equation only depends on genus-0 data.
\end{pro}
{\bf Remark}: Since idempotents span the space of primary vector fields,
combined with the genus-1 topological recursion relation, one can easily
obtain derivatives of any vector field on the big phase space by using the
standard decomposition. Therefore this formula determines the genus-1 generating
function up to an additive constant.

\noindent {\bf Proof of Proposition~\ref{pro:g1ei}}:
It was proved in \cite[Proposition 3.1]{L1} that for any vector fields $v_{i}$,
\begin{eqnarray}
    G_{1}(v_{1}, v_{2}, v_{3}, v_{4})
    &=& \sum_{g \in S_{4}} \left\{
        3 \{v_{g(1)} \bullet v_{g(2)} \} \left<\left<
            \{v_{g(3)} \bullet v_{g(4)} \}
                     \right>\right>_{1}   \right. \nonumber \\
    && \hspace{30pt}
    - 4 v_{g(4)} \left<\left< \{v_{g(1)} \bullet v_{g(2)} \bullet
            v_{g(3)} \}
                     \right>\right>_{1}       \nonumber \\
    && \hspace{30pt}   \left.
     - 6 \left<\left< \left\{[ v_{g(1)} \bullet v_{g(2)}, \,
            v_{g(3)}] \bullet v_{g(4)} \right\}
                     \right>\right>_{1} \right\}. \label{eqn:G1bracket}
\end{eqnarray}
In particular, since idempotents commutes with each other and
$\ve_{i} \qp \ve_{j} = \delta_{ij} \ve_{i}$, we have
\[ G_{1}(\ve_{i}, \ve_{i}, \ve_{i}, \ve_{j}) = - 24 \ve_{j} \gwiione{\ve_{i}} \]
for all $i$ and $j$. So \eqref{eqn:G} implies
\[ \ve_{j} \gwiione{\ve_{i}} = \frac{1}{24} G_{0}(\ve_{i}, \ve_{i}, \ve_{i}, \ve_{j}).\]
Therefore
\[ \bvx \gwiione{\ve_{i}} = \sum_{j=1}^{N} u_{j} \ve_{j} \gwiione{\ve_{i}}
    = \frac{1}{24} G_{0}(\ve_{i}, \ve_{i}, \ve_{i}, \bvx).
     \]

On the other hand, note
that since $\bvl_{0} = - \bvx$, the standard decomposition of $\vl_{0}$ gives
\[ \gwiione{\bvx} = - \gwiione{\vl_{0}} + \gwiione{ T(\tau_{-}(\vl_{0}))}. \]
The genus-1 $L_{0}$ constraint says that $\gwiione{\vl_{0}}$ is a
constant. So by the genus-1 topological recursion relation,
\[ \gwiione{\bvx}
    = \frac{1}{24} \gwii{\tau_{-}(\vl_{0}) \,\, \ga \,\, \gua}
        \, \, \, + \, \, \,  {\rm constant}. \]
By \cite[Lemma 4.5]{L2},
\[ \nabla_{\ve_{i}} \, \, \tau_{-}(\vl_{0})
    = \tau_{-}( \nabla_{\ve_{i}} \vl_{0})
    = \tau_{-}( \vg * \ve_{i})  = 0. \]
Therefore taking derivative on both sides of the above equation
with respect to $\ve_{i}$, we obtain
\[ \ve_{i} \gwiione{\bvx}
    = \frac{1}{24} \gwii{\ve_{i} \,\,
        \tau_{-}(\vl_{0}) \,\, \ga \,\, \gua}. \]
The lemma then follows from the fact that $\gwiione{\ve_{i}} =
\ve_{i} \gwiione{\bvx} -  \bvx \gwiione{\ve_{i}}$ since $\ve_{i}
= [\ve_{i}, \, \bvx] $. $\Box$

To express $\phi_{i}$ in terms of rotation coefficients,
we need the following properties for the vector field $\bvx$.
\begin{lem} \label{lem:derbvx}
\begin{eqnarray*}
&(i)& \bvx \sqrt{g_{i}} =
    \sum_{j=1}^{N} u_{j} \sqrt{g_{j}} \, r_{ij}, \\
&(ii)& \bvx r_{ij} = -r_{ij} \, \, \, {\rm if} \, \, \, i \neq j, \\
&(iii)& \bvx r_{ii} = -r_{ii} + \sum_{j} \frac{r_{ij}}{\sqrt{g_{i}g_{i}}} \,
        \left< \tau_{-}(\vl_{0}),  \, \ve_{j}\right>
        - \frac{1}{g_{i}} <\tau_{-}^{2}(\vl_{0}), \, \ve_{i}> , \\
&(iv)& \nabla_{\bvx} \ve_{i} = \left(
    \sum_{j} u_{j} \sqrt{\frac{g_{j}}{g_{i}}} \, r_{ij} \right) \ve_{i}
    + \sum_{j} (u_{j}-u_{i}) r_{ij}  \sqrt{\frac{g_{i}}{g_{j}}}
    \, \ve_{j}.
\end{eqnarray*}
\end{lem}
The first three properties follow easily from similar properties of $\vx$ and $\vl_{0}$ proved
in Section~\ref{sec:rotcoeff} after using standard decomposition of these vector
fields. The last property follows from Lemma~\ref{lem:rotder} since
$\bvx = \sum_{i} u_{i} \ve_{i}$.
We also observe that by \eqref{eqn:eitoga} and \eqref{eqn:psi-inverse},
for any tensor $Q$,
\begin{equation} \label{eqn:gagua}
 Q(\ga, \gua, \cdots)
    = \sum_{i=1}^{N} \frac{1}{g_{i}} \, Q(\ve_{i}, \ve_{i}, \cdots).
 \end{equation}
in particular
\begin{equation} \label{eqn:Delta}
 \Delta := \ga \qp \gua =  \sum_{i=1}^{N} \frac{1}{g_{i}} \, \ve_{i}.
 \end{equation}

\begin{lem} \label{lem:G0Rot} For any $i$,
\begin{eqnarray*}
G_{0}(\ve_{i}, \ve_{i}, \ve_{i}, \bvx)
&=& \sum_{j}  \sqrt{\frac{g_{i}}{g_{j}}} \, r_{ij} - (r_{ii} + \bvx r_{ii})
   + \sum_{j, k} \left\{ - \sqrt{\frac{g_{j}}{g_{k}}} \, u_{j} r_{ij} r_{ik}
    + \frac{\sqrt{g_{i}g_{j}}}{g_{k}} u_{j} r_{ik} r_{jk} \right\} \\
&&    - 12 \sum_{j} (u_{i} - u_{j}) r_{ij}^{2} .
\end{eqnarray*}
\end{lem}
{\bf Proof}: We write
\[ G_{0}(\ve_{i}, \ve_{i}, \ve_{i}, \bvx) = f_{1} + 3f_{2} -6 f_{3} \]
where
\[ f_{1} = \gwii{\ve_{i} \, \ve_{i} \, \ve_{i} \, \gua}
\gwii{\ga \, \bvx \, \gb \, \gub}
+  \gwii{\ve_{i} \, \ve_{i} \, \ve_{i} \, \bvx \, \gua}
\gwii{\ga \, \gb \, \gub}, \]
\[ f_{2} = \gwii{\ve_{i} \, \ve_{i} \, \bvx \, \gua}
\gwii{\ga \, \ve_{i} \, \gb \, \gub}, \]
and
\[ f_{3} =  \gwii{\ve_{i} \, \ve_{i} \, \gua \, \gub}
\gwii{\ga \, \gb \, \ve_{i} \, \bvx}. \]
We first note that
\[ f_{1} = \bvx \gwii{\ve_{i} \, \ve_{i} \, \ve_{i} \, \Delta}
- 3 \gwii{(\nabla_{\bvx} \ve_{i}) \, \ve_{i} \, \ve_{i} \, \Delta}. \]
Replacing $\Delta$ by $\sum_{j} \frac{\ve_{j}}{g_{j}}$, then computing
$f_{1}$ using Lemma~\ref{lem:4ptRot} and \ref{lem:derbvx}, we obtain
\begin{eqnarray}
f_{1}
&=& \sum_{j}  \sqrt{\frac{g_{i}}{g_{j}}} \, r_{ij} - (r_{ii} + \bvx r_{ii})
   + \sum_{j, k} \left\{ 2 \sqrt{\frac{g_{j}}{g_{k}}} \, u_{j} r_{ij} r_{ik}
    + \frac{\sqrt{g_{i}g_{j}}}{g_{k}} u_{j} r_{ik} r_{jk} \right\}
    \nonumber \\
&&    +3 \sum_{j} (u_{i} - u_{j}) r_{ij}^{2}
    \left( \frac{g_{i}}{g_{j}} - 1 \right).  \label{eqn:f1G0}
\end{eqnarray}
By \eqref{eqn:gagua},
\[ f_{2} = \sum_{j, k} \frac{1}{g_{j} g_{k}}
    \gwii{\ve_{i} \, \ve_{i} \, \bvx \, \ve_{j}}
\gwii{\ve_{j} \, \ve_{i} \, \ve_{k} \, \ve_{k}} \]
and
\[ f_{3} = \sum_{j, k} \frac{1}{g_{j} g_{k}}
    \gwii{\ve_{i} \, \ve_{i} \, \ve_{j} \, \ve_{k}}
\gwii{\ve_{j} \, \ve_{k} \, \ve_{i} \, \bvx}. \]
Since $\bvx = \sum_{k} u_{k} \ve_{k}$,  by Lemma~\ref{lem:4ptRot},
\[ \gwii{\bvx \, \ve_{i} \, \ve_{i} \, \ve_{i}}
     = - \sum_{j} u_{j} r_{ij} \sqrt{g_{i} g_{j}} \]
for any $i$ and
\[ \gwii{\bvx \, \ve_{i} \, \ve_{i} \, \ve_{j}}
     = ( u_{j} - u_{i}) r_{ij} \sqrt{g_{i} g_{j}} \]
if $i \neq j$. Computing $f_{2}$ and $f_{3}$ using these two formulas and
Lemma~\ref{lem:4ptRot}, we obtain
\begin{equation} \label{eqn:f2G0}
 f_{2}
= 2 \sum_{j} \sqrt{\frac{g_{j}}{g_{i}}} \, u_{j} r_{ij} r_{ii}
    - \sum_{j, k} \sqrt{\frac{g_{j}}{g_{k}}} \, u_{j} r_{ij} r_{ik}
    + \sum_{j} (u_{i} - u_{j}) r_{ij}^{2}
    \left( \frac{g_{i}}{g_{j}} + 1 \right)
\end{equation}
and
\begin{equation} \label{eqn:f3G0}
f_{3}
= \sum_{j} \sqrt{\frac{g_{j}}{g_{i}}} \, u_{j} r_{ij} r_{ii}
    + \sum_{j} (u_{i} - u_{j}) r_{ij}^{2}
    \left( \frac{g_{i}}{g_{j}} + 2 \right).
\end{equation}
The lemma then follows from equations (\ref{eqn:f1G0}),
(\ref{eqn:f2G0}), and (\ref{eqn:f3G0}).
$\Box$

\begin{lem} \label{lem:t-L04pt}
\begin{eqnarray*}
\gwii{\ve_{i} \, \tau_{-}(\vl_{0}) \, \ga \, \gua}
&=& - (r_{ii} + \bvx r_{ii})
   + \sum_{j, k} \left\{ - \sqrt{\frac{g_{j}}{g_{k}}} \, u_{j} r_{ij} r_{ik}
    + \frac{\sqrt{g_{i}g_{j}}}{g_{k}} u_{j} r_{ik} r_{jk} \right\}.
\end{eqnarray*}
\end{lem}
{\bf Proof}:
Since $\tau_{-}(\vl_{0}) = \overline{\tau_{-}(\vl_{0})}
    + T(\tau_{-}^{2}(\vl_{0}))$,
the property of $T$ stated in Section~\ref{sec:idem} implies that
\[\gwii{\ve_{i} \, \tau_{-}(\vl_{0}) \, \ga \, \gua}
    = \gwii{\ve_{i} \, \overline{\tau_{-}(\vl_{0})} \, \ga \, \gua}
    + \gwii{\{\ve_{i} \qp \tau_{-}^{2}(\vl_{0})\} \, \ga \, \gua}.
\]
By associativity of the quantum product,
the second term is
\[ \gwii{\{\ve_{i} \qp \tau_{-}^{2}(\vl_{0})\} \, \ga \, \gua}
    = \gwii{\vs \, \tau_{-}^{2}(\vl_{0}) \, \{\ve_{i} \qp \Delta \} }
    = \frac{1}{g_{i}} <\tau_{-}^{2}(\vl_{0}), \, \ve_{i}>.\]
To compute the first term, we note that
$\overline{\tau_{-}(\vl_{0})} =
    \sum_{j=1}^{N}
        \frac{<\overline{\tau_{-}(\vl_{0})}, \, \ve_{j}>}{g_{j}}
            \, \ve_{j}$.
So by Lemma~\ref{lem:4ptRot}, we have
\begin{eqnarray*}
\gwii{\ve_{i} \, \overline{\tau_{-}(\vl_{0})} \, \ga \, \gua}
&=&\sum_{j, k}
        \frac{<\overline{\tau_{-}(\vl_{0})}, \, \ve_{j}>}{g_{j} g_{k}}
        \gwii{\ve_{i} \, \ve_{j} \, \ve_{k} \, \ve_{k}} \\
&=& - <\overline{\tau_{-}(\vl_{0})}, \,
        \sum_{j} \frac{r_{ij}}{\sqrt{g_{i}g_{j}}} \, \ve_{j}>  \\
&&    + \sum_{j} \sqrt{\frac{g_{i}}{g_{j}}} \, r_{ij}
        \left(\frac{<\overline{\tau_{-}(\vl_{0})}, \, \ve_{i}>}{g_{i}}
        - \frac{<\overline{\tau_{-}(\vl_{0})}, \, \ve_{j}>}{g_{j}} \right).
\end{eqnarray*}
By \cite[Corollary 4.9]{L2}, Lemma~\ref{lem:tau-S}, and Corollary~\ref{cor:FiEj},
\[ \overline{\tau_{-}(\vl_{0})} = - \sum_{j} \left( \frac{3}{2} +
    \sum_{i} u_{i} r_{ij} \sqrt{\frac{g_{i}}{g_{j}}} \right) \ve_{j}. \]
Combining above formulas, we have
\begin{eqnarray*}
\gwii{\ve_{i} \, \tau_{-}(\vl_{0}) \, \ga \, \gua}
&=& \frac{1}{g_{i}} <\tau_{-}^{2}(\vl_{0}), \, \ve_{i}>
    - <\overline{\tau_{-}(\vl_{0})}, \,\,
        \sum_{j} \frac{r_{ij}}{\sqrt{g_{i}g_{j}}} \, \ve_{j}>  \\
&&+ \sum_{j, k} \left\{ - \sqrt{\frac{g_{j}}{g_{k}}} \, u_{j} r_{ij} r_{ik}
    + \frac{\sqrt{g_{i}g_{j}}}{g_{k}} u_{j} r_{ik} r_{jk} \right\}.
\end{eqnarray*}
The lemma then follows from Lemma~\ref{lem:derbvx} (iii).
$\Box$

Lemma~\ref{lem:G0Rot} and Lemma~\ref{lem:t-L04pt}
imply the following
\begin{pro} \label{pro:eiF1Rot}
\begin{eqnarray*}
24 \phi_{i} = 24 \gwiione{\ve_{i}}
&=& 12 \sum_{j} (u_{i} - u_{j}) r_{ij}^{2}
    - \sum_{j}  \sqrt{\frac{g_{i}}{g_{j}}} \, r_{ij}.
\end{eqnarray*}
\end{pro}

Let $J = (J_{i \alpha})_{N \times N}$ be the transition matrix from
frame $\{\ga \mid \alpha=1, \ldots, N \}$ to
$\{\ve_{i} \mid i=1, \ldots, N \}$, i.e. $\ve_{i} = J_{i \alpha} \ga$. Then
$J_{i \alpha} = \sqrt{g_{i}} \psi_{i \beta} \eta^{\beta \alpha}$ by
\eqref{eqn:guatoei}. By \eqref{eqn:psipsi-1},
\[ (\det \psi)^{2} = \det \eta \equiv {\rm constant}. \]
So
\[ \det J = \left( \prod_{j=1}^{N} \sqrt{g_{j}} \right) (\det \psi)
        (\det \eta^{-1})
    = ({\rm constant}) \cdot \prod_{j=1}^{N} \sqrt{g_{j}}. \]
Therefore for every $i$,
\[ \ve_{i} \log \det J = \sum_{j=1}^{N} \ve_{i} \log \sqrt{g_{j}}
    = \sum_{j=1}^{N} \frac{\sqrt{g_{i}} \, r_{ij}}{\sqrt{g_{j}}}. \]
So when restricted to the small phase space, the formula in Proposition~\ref{pro:eiF1Rot}
coincides with the formula in \cite{DZ1}.
It is very surprising that $\gwiione{\ve_{i}}$ on the big phase space
has the exactly the same form as on the small phase space. In contrast,
$\gwii{\ve_{i} \, \tau_{-}(\vl_{0}) \, \ga \, \gua}$, which is equal to 0
on the small phase space, and $G_{0}(\ve_{i}, \ve_{i}, \ve_{i}, \bvx)$
have much more complicated expressions on the big phase space.

We can use Proposition~\ref{pro:eiF1Rot} to give a simple proof to the genus-1
Virasoro conjecture for manifolds with semisimple quantum cohomology
(cf. \cite{DZ2} and \cite{L1} for earlier proofs). What is new in this proof
is that we are working directly on the big phase space, while earlier
proofs worked on the small phase space first, then using the constitutive
relation to obtain the result on the big phase space.
Starting from genus-2, the constitutive relation does not exist.
Therefore the approach here will be easier for the generalization to higher
genus cases.
We only give a proof for the genus-1 $L_{1}$-constraint here
 which has the form
\[
 24 \gwiione{\bvx^{2}} = \gwii{\vs \, \vs \, \vv}
\]
where
\[ \vv :=  \tau_{-}(\vl_{1}) \qp \Delta
    + 6 (\vg*(\bvx \qp \gua)) \qp (\vg * \ga)
    + 6 \bvx \qp (\vg* \gua) \qp (\vg * \ga) \]
(cf. \cite[Section 5.1]{L2}). The vector field $\vl_{1}$ is one of the Virasoro vector
fields defined in \cite[Equation (41)]{L2}. The proof for the genus-1 $L_{2}$-constraint
(hence the genus-1 Virasoro conjecture)
can be done similarly. On the other hand, it was proven in \cite{L1} that
the genus-1 $L_{1}$-constraint implies the genus-1 Virasoro conjecture
for all projective varieties.

By \cite[Corollary 4.9]{L2}, Lemma~\ref{lem:tau-S}, and Corollary~\ref{cor:FiEj},
\[ \overline{\tau_{-}(\vl_{1})} = - \sum_{j} \left( 3u_{j} +
    \sum_{i} u_{i}^{2} r_{ij} \sqrt{\frac{g_{i}}{g_{j}}} \right) \ve_{j}. \]
Therefore by \eqref{eqn:gagua} and Corollary~\ref{cor:G*rot},
\[ \vv = \sum_{i, j} \left\{6(u_{i}+u_{j})(u_{i}-u_{j})^{2}r_{ij}^{2}
   - u_{i}^{2} r_{ij} \sqrt{\frac{g_{i}}{g_{j}}} \right\}
   \frac{1}{g_{j}} \ve_{j}. \]
Since $\gwii{\vs \, \vs \, \ve_{j}} = g_{j}$, the prediction of genus-1
$L_{1}$-constraint is
\[
 24 \gwiione{\bvx^{2}} =
 \sum_{i, j} 6(u_{i}+u_{j})(u_{i}-u_{j})^{2}r_{ij}^{2}
   - u_{i}^{2} r_{ij} \sqrt{\frac{g_{i}}{g_{j}}}.
\]
Since $r_{ij} = r_{ji}$ for all $i$ and $j$, we have
\[ \sum_{i, j} 6(u_{i}+u_{j})(u_{i}-u_{j})^{2}r_{ij}^{2}
   =  \sum_{i, j} 12 u_{i}^{2}(u_{i}-u_{j})r_{ij}^{2}. \]
On the other hand, we can compute
$\gwiione{\bvx^{2}} = \sum_{i}u_{i}^{2} \gwiione{\ve_{i}}$
using Proposition~\ref{pro:eiF1Rot}, and the formula obtained is
exactly the same formula predicted by the genus-1 $L_{1}$-constraint.
This finishes the proof of the genus-1 $L_{1}$ constraint for the semisimple
case.

This method also works for studying genus-2 Gromov-Witten invariants.
The idea is similar, though the computation is considerably more complicated.
If the quantum cohomology of the underlying manifold is
semisimple, we can solve the genus-2 universal equations and obtain
the following formula for the genus-2 generating function:
\begin{equation} \label{eqn:F2}
 F_{2} = \frac{1}{2} A_{1}(\tau_{-}(\vs)) +
    \frac{1}{3} A_{1}(\tau_{-}^{2}(\vl_{0}))
    - \frac{1}{6} \sum_{i=1}^{N} u_{i} B(\ve_{i}, \ve_{i}, \ve_{i})
\end{equation}
where $A_{1}$ is the genus-0 and genus-1 part of a genus-2 topological recursion relation of the form
\[ \gwiitwo{T^{2}(\vw)} = A_{1}(\vw) \]
and $B$ is the genus-0 and genus-1 part of an equation due to Belorousski and Pandharipande
(see \cite[Section 2]{L2} for the precise definition of these two tensors).
This formula is much simpler than the one obtained in \cite{L2} which does not use
idempotents. Taking derivatives of $F_{2}$ along $\vl_{1}$, and using rotation coefficients,
we can check that the result coincides with the prediction of the genus-2 $L_{1}$-constraint.
As proved in \cite{L2}, this
implies that the genus-2 Virasoro conjecture is true for manifolds with semisimple
quantum cohomology. The details will be given in a separate paper \cite{L4}.

It is also interesting to compare the solutions to universal equations at genus 1
and genus 2. For this purpose, we need the following genus-0 relation:
\begin{lem} \label{lem:g0g1}
\begin{eqnarray*}
&(i) & G_{0}(\vw, \ve_{i}, \ve_{j}, \ve_{k})=0 \hspace{20pt}
        {\rm for \,\,\, any \,\,\, vector \,\,\, field} \,\,\, \vw
        \,\,\, {\rm if} \,\,\, i \neq j \neq k.  \\
& (ii) & G_{0}(\ve_{i}, \ve_{i}, \ve_{i}, \ve_{j})
        = G_{0}(\ve_{i}, \ve_{j}, \ve_{j}, \ve_{j})
        = - G_{0}(\ve_{i}, \ve_{i}, \ve_{j}, \ve_{j}) \hspace{20pt}
            {\rm if} \,\,\, i \neq j.
\end{eqnarray*}
\end{lem}
{\bf Proof}: By \eqref{eqn:G}, the formulas in this lemma are
equivalent to the corresponding formulas with $G_{0}$ replaced by
$G_{1}$. For $i \neq j \neq k$,
\begin{eqnarray*}
 G_{1}(\vw, \ve_{i}, \ve_{j}, \ve_{k})
   & = &
   -4\gwii{(\vw \qp \ve_{i}) \, \ve_{j} \, \ve_{k} \, \gua} \gwiione{\ga}
    + 12\gwii{\vw \, \ve_{j} \, \ve_{k} \, \gua}
            \gwiione{(\ga \qp \ve_{i})}   \\
 && + {\rm cyclic \,\,\, permutation \,\,\, of \,\,\, these \,\,\,
        two \,\,\, terms \,\,\, with \,\,\, respect \,\,\, to \,\,\,}
    i, j, k.
\end{eqnarray*}
Note that $\vw \qp \ve_{i}$ is proportional to $\ve_{i}$.
Lemma~\ref{lem:derejfj} (c) and \eqref{eqn:weiej4pt} imply that
$G_{1}(\vw, \ve_{i}, \ve_{j}, \ve_{k}) =0$. This proves (i).

By \eqref{eqn:G1bracket},
\[ G_{1}(\ve_{i}, \ve_{i}, \ve_{i}, \ve_{j})
    = - 24 \ve_{j} \gwiione{\ve_{i}}  \]
and
\[ G_{1}(\ve_{i}, \ve_{i}, \ve_{j}, \ve_{j})
    = 12 \ve_{j} \gwiione{\ve_{i}} + 12 \ve_{i} \gwiione{\ve_{j}}. \]
Since $[ \ve_{i}, \, \ve_{j}]=0$, $\ve_{j} \gwiione{\ve_{i}} =
\ve_{i} \gwiione{\ve_{j}}$. Together with the fact that $G_{1}$
and $G_{0}$ are symmetric tensors, this proves (ii). $\Box$

{\bf Remark}: Comparing to Lemma~\ref{lem:derejfj} (c) and
Lemma~\ref{lem:4ptg0} (c), we see that when applied to
idempotents, $G_{0}$ and $G_{1}$ has very similar properties as
the genus-0 4-point functions.

Lemma~\ref{lem:g0g1} (ii)
implies that
\[ G_{0}(\ve_{i}, \ve_{i}, \ve_{i}, \bvx) =
    \sum_{j} u_{j} G_{0}(\ve_{i}, \ve_{i}, \ve_{i}, \ve_{j})
    = \sum_{j} u_{j} G_{0}(\ve_{j}, \ve_{j}, \ve_{j}, \ve_{i}). \]
Hence Proposition~\ref{pro:g1ei} can be rewritten as
\[ \gwiione{\ve_{i}}
    = \frac{1}{24} \left( \gwii{\ve_{i} \, \tau_{-}(\vl_{0}) \, \ga \, \gua}
        - \sum_{j} u_{j} G_{0}(\ve_{j}, \ve_{j}, \ve_{j}, \ve_{i}) \right).\]
If we ignore $\ve_{i}$ on both sides  and notice that the first
term on the right hand side of this equation comes from the
genus-0 part of the genus-1 topological recursion relation, this
formula is very similar to the expression for $F_{2}$  in
\eqref{eqn:F2}. Based on this observation, we might
speculate that higher genus generating functions $F_{g}$ with $g>2$
can be solved in a similar fashion.


\vspace{30pt}
\noindent
Department of Mathematics  \\
University of Notre Dame \\
Notre Dame,  IN  46556, USA \\

\vspace{10pt} \noindent E-mail address: {\it xliu3@nd.edu}

\end{document}